\documentclass[11pt,a4paper]{amsart}

\usepackage{amssymb,amsmath,graphics,verbatim}
\usepackage{latexsym}
\usepackage{eucal}
\usepackage{a4wide}

\newtheorem{proposition}{Proposition}[section]
\newtheorem{theorem}{Theorem}[section]
\newtheorem{lemma}[theorem]{Lemma}

\newtheorem{coro}[theorem]{Corollary}

\newcommand{\mc}{\mathcal}
\newcommand{\cal}{\mathcal}
\newcommand{\rr}{\mathbb{R}}
\newcommand{\R}{\mathbb{R}}
\newcommand{\nn}{\mathbb{N}}
\newcommand{\cc}{\mathbb{C}}

\newcommand{\zz}{\mathbb{Z}}

\newcommand{\eps}{\epsilon}

\newcommand{\pl}{\partial}
\newcommand{\x}{\times}

\newcommand{\bbar}{\overline}

\newcommand{\cjd}{\rangle}
\newcommand{\cjg}{\langle}

\newcommand{\demi}{\frac{1}{2}}

\def\qed{\hfill$\square$\medskip}

\begin{document}
\title[Calder\'on Problem on Riemann surfaces]{Calder\'on inverse Problem for the Schr\"odinger Operator on Riemann Surfaces}
\author{Colin Guillarmou}
\address{Laboratoire J.A. Dieudonn\'e\\
U.M.R. 6621 CNRS\\
Universit\'e de Nice Sophia-Antipolis\\
Parc Valrose, 06108 Nice\\France}
\email{cguillar@math.unice.fr}

\author{Leo Tzou}
\address{Department of Mathematics\\
Stanford University\\
Stanford, CA 94305, USA.}
\email{ltzou@math.stanford.edu.}

\begin{abstract}
On a fixed smooth compact Riemann surface with boundary $(M_0,g)$, we show that the Cauchy data space (or Dirichlet-to-Neumann map $\mc{N}$) 
of the Schr\"odinger operator $\Delta +V$ with $V\in C^2(M_0)$ 
determines uniquely the potential $V$. We also discuss briefly the corresponding consequences for potential scattering at $0$ frequency
on Riemann surfaces with asymptotically Euclidean or asymptotically hyperbolic ends.
\end{abstract}

\maketitle

\begin{section}{Introduction}
The problem of determining the potential in the Schr\"odinger operator by boundary measurement  
goes back to Calder\'on \cite{Ca}. Mathematically, it amounts to ask if one can detect some data 
from boundary measurement in a domain (or manifold) $\Omega$ with boundary. The typical model to have in mind is the Schr\"odinger operator
$P=\Delta_g+V$ where $g$ is a metric and $V$ a potential, then we define the Cauchy data space by 
\[\mc{C}:=\{(u|_{\pl\Omega},\pl_nu|_{\pl\Omega})\in C^\infty(\pl\Omega)\x C^\infty(\pl\Omega)\,; \,u \in \ker P\}\]
where $\pl_n$ is the interior pointing normal vector field to $\pl\Omega$.\\ 

The first natural question is the following \emph{full data} inverse problem: does the Cauchy data space determine uniquely the metric
$g$ and/or the potential $V$? 
In a sense, the most satisfying known results are when the domain $\Omega\subset \rr^n$ is already known and $g$ is
the Euclidean metric, then the recovery of $V$ has been proved in dimension $n>2$ by Sylvester-Uhlmann \cite{SU}
and very recently in dimension $2$ by Bukgheim \cite{Bu}. A related question is the conductivity problem 
which consists in taking $V=0$ and replacing $\Delta_g$ by $-{\rm div}\sigma\nabla$ where $\sigma$
is a definite positive symmetric matrix. An elementary observation shows that 
the problem of recovering an sufficiently smooth isotropic conductivity (i.e. $\sigma=\sigma_0{\rm Id}$ for a function $\sigma_0$) 
is contained in the problem above of recovering a potential $V$. 
For domain of $\rr^2$, Nachman \cite{N} used the $\bar{\pl}$ techniques to show 
that the Cauchy data space determines the conductivity. Recently a new approach developed by Astala and P\"aiv\"arinta in \cite{AP} improved this result to assuming that the conductivity is only a $L^\infty$ scalar function. This was later generalized to $L^\infty$ anisotropic conductivities by Astala-Lassas-P\"aiv\"arinta in \cite{ALP}. We notice that there still 
are rather few results in the direction of recovering the Riemannian manifold $(\Omega,g)$ 
when $V=0$, for instance the surface case by Lassas-Uhlmann \cite{LU} (see also \cite{Be,HM}), the real-analytic manifold case 
by Lassas-Taylor-Uhlmann \cite{LTU} (see also \cite{GSB} for the Einstein case),  the case of manifolds admitting limiting Carleman weights 
and in a same conformal class by Dos Santos Ferreira-Kenig-Salo-Uhlmann \cite{DKSU}.\\

The second natural, but harder, problem is the \emph{partial data} inverse problem: if $\Gamma_1$ and $\Gamma_2$ are open subsets of $\pl\Omega$, 
does the partial Cauchy data space 
\[\mc{C}_{\Gamma_{1},\Gamma_2}:=\{(u_{\pl\Omega},\pl_nu|_{\Gamma_1})\in C^\infty(\partial\Omega)\x C^\infty(\Gamma_1)\,; \,u \in \ker P;\, u|_{\pl\Omega}\in C^{\infty}_0(\Gamma_2)\}\]
determine the domain $\Omega$, the metric, the potential?
For a fixed domain of $\rr^n$, the recovery of the potential if $n>2$ with partial data measurements was initiated by Bukhgeim-Uhlmann \cite{BU} and later improved by 
Kenig-Sj\"ostrand-Uhlmann \cite{KSU} to the case where $\Gamma_1$ and $\Gamma_2$ are respectively open subsets of the "front" and "back" ends of the domain. We refer the reader to the references for a more precise formulation of the problem. In dimension $2$, the recent works of Imanuvilov-Uhlmann-Yamamoto \cite{IUY,IUY2} solves the
problem for fixed domains of $\rr^2$ in the case when $\Gamma_1 = \Gamma_2$ and $\Gamma_1 = {\overline\Gamma_2}^c$.\\

In this work, we address the same question when the background domain is a known Riemann surface with boundary. 
We prove the following recovery result under full data measurement (the partial data case will be included in a forthcoming paper): 
\begin{theorem}
\label{identif}
Let $(M_0,g)$ be a smooth compact Riemann surface with boundary and let $\Delta_g$ be its positive Laplacian. 
Let $V_1,V_2\in C^2(M_0)$ be two real potentials and let $\mc{C}_1,\mc{C}_2$ be the respective Cauchy data spaces.
If $\mc{C}_1=\mc{C}_2$ then $V_1=V_2$.
\end{theorem}
\end{section}

Notice that when $\Delta_g+V_i$ do not have $L^2$ eigenvalues for the Dirichlet condition,
the statement above can be given in terms of Dirichlet-to-Neumann operators. 
Since $\Delta_{\hat{g}}=e^{-2\varphi}\Delta_g$ when $\hat{g}=e^{2\varphi}g$ for some function $\varphi$, it is clear that 
in the statement in Theorem \ref{identif}, we only need to fix the conformal class of $g$ instead of the metric $g$ 
(or equivalently to fix the complex structure on $M$). Observe also that Theorem \ref{identif} implies that, for a fixed Riemann surface with boundary $(M_0,g)$, the Dirichlet-to-Neumann map for the operator $u\to -{\rm div}_{g}(\gamma \nabla^gu)$ determines the isotropic conductivity $\gamma$ if $\gamma\in C^{4}(M_0)$ in the sense that two conductivities giving rise to the same Dirichlet-to-Neumann are equal. 
This is a standard observation by transforming the conductivity problem to a potential problem 
with potential $V:=(\Delta_g \gamma^{\demi})/\gamma^\demi$. So our result also extends that of Henkin-Michel \cite{HM} in the case
of isotropic conductivities.\\ 

The method to reconstruct the potential follows \cite{Bu,IUY} and is based on the construction of special complex geometric optic solutions of 
$(\Delta_g+V)u=0$, more precisely solutions of the form $u=e^{\Phi/h}(a+r(h))$ where $h>0$ is a small parameter, $\Phi$ and $a$ are holomorphic functions on $(M,g)$ and $r(h)$ is an error term small as $h\to 0$. The idea of \cite{Bu} to reconstruct $V(p)$ for $p\in M$
is to take $\Phi$ with a non-degenerate critical point at $p$ and then use staionary phase as $h\to 0$. In our setting,
one of our main contribution is the construction of the holomorphic Carleman weights $\Phi$ which is quite more complicated 
since we are working on a Riemann surface instead of a domain of $\cc$. We also need to prove a Carleman estimate on the surface
for this degenerate weight.\\

A trivial consequence we mention is an inverse result for the scattering operator on asymptotically hyperbolic (AH in short) surfaces for potential decaying at the boundary. Recall that an AH surface is an open complete Riemannian surface $(X,g)$ such that $X$ is the interior of a smooth compact surface with boundary $\bar{X}$, and for any smooth boundary defining function $x$ of $\pl\bar{X}$, $\bar{g}:=x^2g$
extends as a smooth metric to $\bar{X}$, with curvature tending to $-1$ at $\pl\bar{X}$.
If $V\in C^{\infty}(\bar{X})$ and $V=O(x^2)$, then we can define a scattering map as follows (see for instance \cite{JSB,GRZ} or \cite{GG}):
first the $L^2$ kernel $\ker_{L^2}(\Delta_g+V)$ is a finite dimensional subspace of $xC^{\infty}(\bar{X})$ and 
in one-to-one correspondence with $E:=\{(\pl_x\psi)|_{\pl\bar{X}}; \psi\in\ker_{L^2}(\Delta_g+V)\}$ where $\pl_x:=\nabla^{\bar{g}} x$ is the normal vector field to $\pl\bar{X}$ for $\bar{g}$, then  
for $f\in C^{\infty}(\pl\bar{X})/E$, there exists a unique function $u\in C^{\infty}(\bar{X})$
such that $(\Delta_g+V)u=0$ and $u|_{\pl\bar{X}}=f$, then one can see that the scattering map 
$\mc{S}:C^\infty(\pl\bar{X})/E\to C^{\infty}(\pl\bar{X})/E$ is defined by $\mc{S}f:=\pl_x u|_{\pl\bar{X}}$. We thus obtain 
\begin{coro}\label{coroAH}
Let $(X,g)$ be an asymptotically hyperbolic manifold and $V_1,V_2\in x^{2}C^{\infty}(\bar{X})$ two potentials such that 
$\Delta_g+V_i$. Assume that 
\[\{\pl_xu|_{\pl\bar{X}}; u\in \ker_{L^2}(\Delta_g+V_1)\}
=\{\pl_xu|_{\pl\bar{X}}; u\in \ker_{L^2}(\Delta_g+V_2)\}\] 
and let $\mc{S}_j$ be the scattering map for the operator $\Delta_g+V_j$, $j=1,2$, then $\mc{S}_1=\mc{S}_2$ implies that $V_1=V_2$. 
\end{coro}

To conclude, we consider the case of inverse scattering at $0$ frequency 
for $\Delta+V$ with $V$ compactly supported on an asymptotically Euclidean surface.
An asymptotically Euclidean surface is a non-compact Riemann surface $(X,g)$, which compactifies into $\bar{X}$ and such that 
the metric in a collar $(0,\eps)_x\x \pl\bar{X}$ near the boundary is of the form  
\[g=\frac{dx^2}{x^4}+\frac{h(x)}{x^2}\]
where $h(x)$ is a smooth one-parameter family of metrics on $\pl\bar{X}$ with $h(0)=d\theta^2_{S^1}$ is the metric 
with length $2\pi$ on each copy of $S^1$ that forms the connected components of $\pl\bar{X}$. 
Notice that using the coordinates $r:=1/x$, $g$ is asymptotic to $dr^2+r^2d\theta^2_{S^1}$ near $r\to \infty$.
A particular case is given by the surfaces with Euclidean ends, 
ie. ends isometric to $\rr^2\setminus B(0,R)$ where $B(0,R)=\{z\in\rr^2; |z|\geq R\}$.
The elements of kernel of $\Delta_g+V$ which have polynomial growth of order, say, $N\in \rr^+$
at $\infty$ form a finite dimensional vector space of dimension increasing with $N$, 
so it is clear that one can not expect to recover much 
information from this space. Thus in a sense, 
the Cauchy data space at $0$ frequency can be interpreted as the coefficients
in the asymptotic expansion of tempered solutions of $(\Delta_g+V)u=0$, we mean the 
solutions which lie in a weighted space $x^{\alpha}L^2$ for some $\alpha\in\rr$.
We show 
\begin{theorem}\label{asconic}
Let $(X,g)$ be an asymptotically Euclidean surface and let $V_1,V_2$ be two compactly supported 
$C^2$ potentials in $X$. Assume that for all $\alpha\in \R$ and for any element $\psi$ in $\ker_{x^{\alpha}L^2}(\Delta_g+V_1)$ 
there is a $\varphi\in\ker_{x^{\alpha}L^2}(\Delta_g+V_2)$ such that $\varphi-\psi=O(x^\infty)$
and conversely. Then $V_1=V_2$.
\end{theorem}
We remark that a similar statement could be proved along the same lines without difficulties for asymptotically 
conic manifolds but we prefer to restrict to the asymptotically Euclidean case for simplicity.

Another straightforward corollary in the asymptotically Euclidean case 
is the recovery of a compactly supported potential from the scattering operator at a positive frequency. The proof
is essentially the same as for the operator $\Delta_{\rr^n}+V$ once one know Theorem \eqref{identif}. 

\begin{section}{Harmonic and Holomorphic Morse Functions on a Riemann Surface}

\subsection{Riemann surfaces}
Let $(M_0,g_0)$ be a compact connected smooth Riemannian surface with boundary $\pl M_0$. The surface $M_0$ can be considered as 
a subset of a compact Riemannian surface $(M,g)$, for instance by taking the double of $M_0$ and extending smoothly 
the metric $g_0$ to $M$. 
The conformal class of $g$ on the closed surface $M$ induces a structure of closed Riemann surface, i.e. a closed surface equipped with a complex structure via holomorphic charts $z_\alpha:U_{\alpha}\to \cc$. The Hodge star operator $\star$ acts on the cotangent bundle $T^*M$, its eigenvalues are 
$\pm i$ and the respective eigenspace $T_{1,0}^*M:=\ker (\star+i{\rm Id})$ and $T_{0,1}^*M:=\ker(\star -i{\rm Id})$
are sub-bundle of the complexified cotangent bundle $\cc T^*M$ and the splitting $\cc T^*M=T^*_{1,0}M\oplus T_{0,1}^*M$ holds as complex vector spaces.
Since $\star$ is conformally invariant on $1$-forms on $M$, the complex structure depends only on the conformal class of $g$.
In holomorphic coordinates $z=x+iy$ in a chart $U_\alpha$,
one has $\star(udx+vdy)=-vdx+udy$ and
\[T_{1,0}^*M|_{U_\alpha}\simeq \cc dz ,\quad T_{0,1}^*M|_{U_\alpha}\simeq \cc d\bar{z}  \]
where $dz=dx+idy$ and $d\bar{z}=dx-idy$. We define the natural projections induced by the splitting of $\cc T^*M$ 
\[\pi_{1,0}:\cc T^*M\to T_{1,0}^*M ,\quad \pi_{0,1}: \cc T^*M\to T_{0,1}^*M.\]
The exterior derivative $d$ defines the De Rham complex $0\to \Lambda^0\to\Lambda^1\to \Lambda^2\to 0$ where $\Lambda^k:=\Lambda^kT^*M$
denotes the real bundle of $k$-forms on $M$. Let us denote $\cc\Lambda^k$ the complexification of $\Lambda^k$, then
the $\pl$ and $\bar{\pl}$ operators can be defined as differential operators 
$\pl: \cc\Lambda^0\to T^*_{1,0}M$ and $\bar{\pl}:\cc\Lambda_0\to T_{0,1}^*M$ by 
\[\pl f:= \pi_{1,0}df ,\quad \bar{\pl}:=\pi_{0,1}df,\]
they satisfy $d=\pl+\bar{\pl}$ and are expressed in holomorphic coordinates by
\[\pl f=\pl_zf\, dz ,\quad \bar{\pl}f=\pl_{\bar{z}}f \, d\bar{z}.\]  
with $\pl_z:=\demi(\pl_x-i\pl_y)$ and $\pl_{\bar{z}}:=\demi(\pl_x+i\pl_y)$.
Similarly, one can define the $\pl$ and $\bar{\pl}$ operators from $\cc \Lambda^1$ to $\cc \Lambda^2$ by setting 
\[\pl (\omega_{1,0}+\omega_{0,1}):= d\omega_{0,1}, \quad \bar{\pl}(\omega_{1,0}+\omega_{0,1}):=d\omega_{1,0}\]
if $\omega_{0,1}\in T_{0,1}^*M$ and $\omega_{1,0}\in T_{1,0}^*M$.
In coordinates this is simply
\[\pl(udz+vd\bar{z})=\pl v\wedge d\bar{z},\quad \bar{\pl}(udz+vd\bar{z})=\bar{\pl}u\wedge d{z}.\]
There is a natural operator, the Laplacian acting on functions and defined by 
\[\Delta f:= -2i\star \bar{\pl}\pl f =d^*d \]
where $d^*$ is the adjoint of $d$ through the metric $g$ and $\star$ is the Hodge star operator mapping 
$\Lambda^2$ to $\Lambda^0$ and induced by $g$ as well.
 
To construct Carleman weights, we will use strongly the Riemann-Roch theorem, so  
for the convenience of the reader we recall it (see Farkas-Kra \cite{FK} for more details). A divisor $D$ on $M$ is an element 
\[D=\big((p_1,n_1), \dots, (p_k,n_k)\big)\in (M\x\zz)^k, \textrm{ where }k\in\nn\]  
which will also be denoted $D=\prod_{i=1}^kp_i^{n_i}$ or $D=\prod_{p\in M}p^{\alpha(p)}$ where $\alpha(p)=0$
for all $p$ except $\alpha(p_i)=n_i$. The inverse divisor of $D$ is defined to be 
$D^{-1}:=\prod_{p\in M}p^{-\alpha(p)}$ and 
the degree of the divisor $D$ is defined by $\deg(D):=\sum_{i=1}^kn_i=\sum_{p\in M}\alpha(p)$. 
A meromorphic function on $M$ is said to have divisor $D$ if $(f):=\prod_{p\in M}p^{{\rm ord}(p)}$ is equal to $D$,
where ${\rm ord}(p)$ denotes the order of $p$ as a pole or zero of $f$ (with positive sign convention for zeros). Notice that 
in this case we have $\deg(f)=0$.
For divisors $D'=\prod_{p\in M}p^{\alpha'(p)}$ and $D=\prod_{p\in M}p^{\alpha(p)}$, 
we say that $D'\geq D$ if $\alpha'(p)\geq \alpha(p)$ for all $p\in M$.
The same exact notions apply for meromorphic $1$-forms on $M$. Then we define for a divisor $D$
\[\begin{gathered}
r(D):=\dim \{f \textrm{ meromorphic functions on } M; (f)\geq D\},\\
i(D):=\dim\{u\textrm{ meromorphic 1 forms on }M; (u)\geq D\}.
\end{gathered}\]
The Riemann-Roch theorem states the following identity: for all divisor $D$ on the closed Riemann surface $M$ of genus $g$, 
\begin{equation}\label{riemannroch}
r(D^{-1})=i(D)+\deg(D)-g+1.
\end{equation}
Notice also that for any divisor $D$ with $\deg(D)>0$, one has $r(D)=0$ since $\deg(f)=0$ for all $f$ meromorphic. 
By \cite[Th. p70]{FK}, let $D$ be a divisor, then for any non-zero meromorphic 1-form $\omega$ on $M$, one has 
\begin{equation}\label{abelian}
i(D)=r(D(\omega)^{-1})
\end{equation}
which is thus independent of $\omega$.
 
\subsection{Morse holomorphic functions with prescribed critical points} 
The main result of this section is the following 
\begin{proposition}
\label{criticalpoints}
Let $q$ be a point in $M\setminus M_0$ and let $\mc{O} \subset M\backslash\{q\}$ be an open subset with smooth boundary of the punctured 
Riemann surface $M\setminus\{q\}$ such that $M_0\subset\mc{O}$. 
Then there exists a dense set of points $p$ in $\mc{O}$ such that there exists a Morse holomorphic function $f$ on $\mc{O}$ 
which has a critical point at $p$.
\end{proposition} 
We first prove an auxiliary result which states that for any point $p\in M\backslash\{q\}$ one can find a holomorphic function 
on $M\backslash\{p,q\}$, meromorphic on $M$ and with a pole or zero of any desired order at $p$.
\begin{lemma} \label{poleszeros}
Let $p\in M\backslash\{q\}$ and $n\in \nn$. Then there exist meromorphic functions $h_n,k_n$ on $M$
such that $k_n$ is holomorphic on $M\backslash\{q,p\}$ with a pole of order $n$ at $p$ and $h_n$ is holomorphic
on $M\setminus\{q\}$ with a zero of order $n$ at $p$.  
\end{lemma}
\noindent{\bf Proof}. First we claim that there exists $N_0\in\nn$ so that for all $N\geq N_0$, there is meromorphic function on $M$, holomorphic on $M\setminus\{p\}$ with a pole of order $N$ at $p$. Indeed, fix a meromorphic 1-form $\omega$, then by \eqref{riemannroch}, we know that
for $D:=p^{N}$ with $N>g-1$, then $r(D^{-1})>0$. Moreover, if $\deg(D(\omega^{-1}))>0$, one has $r(D(\omega)^{-1})=0$ 
so we conclude by \eqref{abelian} and \eqref{riemannroch} that if $N_0>g-1$ is taken large enough and $N\geq N_0$ then 
$r(D^{-1})=N-g+1$, which implies that there is a meromorphic function $f_N$ with a pole of order $N$ at $p$ and no other poles.
By Riemann-Roch again \eqref{riemannroch}, one has $r(D^{-1})>0$ if $D=p^{-\ell}q^{k}$ with $k,l\in\nn$ and $k-\ell>g-1$. Thus
there exists a meromorphic function $h$ on $M$, holomorphic on $M\setminus\{q\}$, 
with a pole of order, say $k'\leq k$, at $q$ and a zero of order, say $\ell'\geq \ell$ at $p$.  
By possibly taking powers of $h$, we may assume that $p$ is a zero of $h$ of order say $N$ with $N>N_0$. 
Then the function $h_n:= (f_{N-1}h)^n$ is meromorphic on $M$, holomorphic on $M\setminus\{q\}$, and with a zero of order $n$ at $p$. 
Similarly, the function $k_n:=(h_{N-1}f_N)^{n}$ is meromorphic on $M$, holomorphic on $M\setminus\{p,q\}$ and with a pole of 
order $n$ at $p$.\qed\\

Fix $k>2$ a large integer, we denote by $C^k(\bar{\mc{O}})$ the Banach space of $C^k$ real valued functions on $\bar{\mc{O}}$. 
Then the set of harmonic functions on $\mc{O}$ which are in the Banach space
$C^{k}(\bar{\mc{O}})$ (and smooth in $\mc{O}$ by elliptic regularity) 
is the kernel of the continuous map $\Delta:C^k(\bar{\mc{O}})\to C^{k-2}(\bar{\mc{O}})$, 
and so it is a Banach subspace of $C^k(\bar{\mc{O}})$. 
%More precisely, there is a continuous injective linear map (called Poisson operator) 
%\[P :H^{k-\demi}(\pl \bar{\mc{O}})\to H^{k}(\bar{\mc{O}})\]
%such that $\Delta Pf=0$ and $Pf|_{\pl\bar{\mc{O}}}=f$. The image of $H^{k-\demi}(\pl\bar{\mc{O}})$ under $P$
%is the set of harmonic functions in $H^k(\bar{\mc{O}})$, this is a Hilbert subspace of $H^k(\bar{\mc{O}})$. 
The set $\mc{H}\subset C^k(\bar{\mc{O}})$ of harmonic function $u$ in $C^k(\bar{\mc{O}})$ such there exists 
$v\in C^{k}(\bar{\mc{O}})$ harmonic with $u+iv$ holomorphic on $\mc{O}$ is a Banach subspace of $C^{k}(\bar{\mc{O}})$ of 
finite codimension. Indeed, let $\{\gamma_1,..,\gamma_N\}$ be a homology basis for $\mc{O}$, then
\[\mc{H}=\ker L , \textrm{ with } L: \ker\Delta\cap C^k(\bar{\mc{O}})\to \cc^N \textrm{ defined by }
L(u):=\Big(\frac{1}{\pi i}\int_{\gamma_j}\pl u\Big)_{j=1,\dots,N}.\]
%Now let $\mc{B}\subset H^{k-\demi}(\pl\bar{\mc{O}})$ be the set of boundary values of functions in $\mc{H}$, this is given by 
%$\mc{B}:=P^{-1}\mc{H}$ where $P^{-1}$ is the inverse of the continuous linear isomorphism 
%$P:H^{k-\demi}(\pl\bar{\mc{O}})\to P(H^{k-\demi}(\pl\bar{\mc{O}}))$.  In other words, 
%we have the isomorphism of Hilbert spaces $P:\mc{B}\to \mc{H}$.
We now show
\begin{lemma}\label{morsedense}
The set of functions $u\in\mc{H}$ which are Morse in $\mc{O}$ 
is dense in $\mc{H}$ with respect to the $C^k(\bar{\mc{O}})$ topology.
\end{lemma}
\noindent{\bf Proof}. We use an argument very similar to those used by Uhlenbeck \cite{Uh}.
We start by defining $m: \bar{\mc{O}}\times \mc{H}\to T^*\mc{O}$ by $(p,u) \mapsto (p,du(p))\in T_p^*\mc{O}$. 
This is clearly a smooth map, linear in the second variable, moreover $m_u:=m(.,u)=(\cdot, du(\cdot))$ is  
Fredholm since $\mc{O}$ is finite dimensional. The map $u$ is a Morse functions if and only if 
$m_u$ is transverse to the zero section, denoted $T_0^*\mc{O}$, of $T^*\mc{O}$, ie. if 
\[\textrm{Image}(D_{p}m_u)+T_{m_u(p)}(T_0^*\mc{O})=T_{m_u(p)}(T^*\mc{O}),\quad \forall p\in \mc{O} \textrm{ such that }m_u(p)=(p,0).\]
which is equivalent to the fact that the Hessian of $u$ at critical points is 
non-degenerate (see for instance Lemma 2.8 of \cite{Uh}). 
We recall the following transversality theorem (\cite[Th.2]{Uh} or \cite{Ab,Qu}):
\begin{theorem}\label{transversality}
Let $m : X\times \mc{H} \to W$ be a $C^k$ map, where $X$, $\mc{H}$, and $W$ are separable Banach manifolds 
with $W$ and $X$ of finite dimension. Let $W'\subset W$ be a submanifold such that $k>\max(1,\dim X-\dim V+\dim V')$.
If $m$ is transverse to $W'$  then the set 
$\{u\in \mc{H}; m_u \textrm{ is transverse to } W'\}$ is dense in $\mc{H}$, more precisely it is a set of second category. 
\end{theorem}
We want to apply it with $X:=\mc{O}$, $W:=T^*\mc{O}$ and $W':=T^*_0\mc{O}$, and the map $m$ is defined above. 
We have thus proved our Lemma if one can show that $m$ is transverse to $W'$. 
Let $(p,u)$ such that $m(p,u)=(p,0)\in V'$. Then identifying $T_{(p,0)}(T^*\mc{O})$ with $T_p\mc{O}\oplus T^*_p\mc{O}$, one has
\[Dm_{(p,u)}(z,v)=(z,dv(p)+{\rm Hess}_p(u)z)\]
where ${\rm Hess}_pu$ is the Hessian of $u$ at the point $p$, viewed as a linear map from $T_p\mc{O}$ to $T^*_p\mc{O}$. 
To prove that $m$ is transverse to $W'$
we need to show that $(z,v)\to (z, dv(p)+{\rm Hess}_p(u)z)$ is onto from $T_p\mc{O}\oplus \mc{H}$ to $T_p\mc{O}\oplus T^*_p\mc{O}$, which 
is realized for instance if the map $v\to dv(p)$ from $\mc{H}$ to  $T_p^*\mc{O}$ is onto.
But from Lemma \ref{poleszeros}, we know that there exists a holomorphic function $v$ on $M\setminus\{q\}$ (thus on $\bar{\mc{O}}$)
such that $v(p)=0$ and $dv(p)\not=0$ as a linear map $T_p\mc{O}\to \cc$, we can then take its real and imaginary parts
$v_1$ and $v_2$, both are real valued harmonic smooth function on $\bar{\mc{O}}$ thus in $\mc{H}$, and $dv_1(p)$ and $dv_2(p)$ are
linearly independent in $T^*_p\mc{O}$ by the Cauchy-Riemann equation $\bar{\pl} v=0$. This shows our claim and ends the proof
by using Theorem \ref{transversality}.\qed\\

\noindent{\bf Proof of Proposition \ref{criticalpoints}} Let $p$ be a point of $\mc{O}$ and let $u$ 
be a holomorphic function with a nondegenerate critical point at $p$, the existence is insured by Lemma \ref{poleszeros}. 
By Lemma \ref{morsedense}, there exist Morse holomorphic functions $(u_j)_{j\in\nn}$ 
such that $u_j\to u$ in $C^k(\bar{\mc{O}},\cc)$ for any fixed $k$ large. 
Let $\eps>0$ small and let $U\subset\mc{O}$ be a neighbourhood containing $p$ and no other critical points of $u$,
and with boundary a smooth circle of radius $\eps$. In complex local coordinates near $p$, 
we can consider $\pl u$ and $\pl u_j$ as holomorphic functions on an open set of $\cc$.
Then by Rouche's theorem, it is clear that $\pl{u_j}$ has precisely one zero in $U$ if $j$ is large enough. 
This completes the proof.\qed\\

\end{section}

\begin{section}{Carleman Estimate for Harmonic Weights with Critical Points}

In this section, we prove a Carleman estimate using harmonic weight with non-degenerate critical points, in 
way similar to \cite{IUY}: 

%\begin{proposition}
%\label{carlemanestimate}
%Let $q\in L^\infty(\cal{M})$ and $\phi$ be a harmonic Morse function. For all $h>0$ small enough, we have that for all $u\in C^{\infty}_0({\cal M})$ the estimate
%\[\|e^{-\phi/h}(\Delta +V)e^{\phi/h}u\|^2 \geq C(\frac{1}{h}\|u\|^2 + \frac{1}{h^2}\|ud\phi\|^2 + \|du\|^2)\]
%\end{proposition}
\begin{proposition}
\label{carlemanestimate}
Let $(M,g)$ be a Riemann surface with boundary, with $\bar{M}:=M\cup\pl\bar{M}$, and let $\varphi: \bar{M}\to\R$ be a harmonic function with non-degenerate critical points. 
Then for all $V\in L^\infty$ there exists an $h_0>0$ such that for all $h<h_0$ and $u\in C^\infty_0(M)$, we have 
\begin{eqnarray}
\label{carleman}
\frac{1}{h}\|u\|^2 + \frac{1}{h^2}\|u |d\varphi|\|^2 + \|du\|^2 \leq C\|e^{-\varphi/h}(\Delta_g + V)e^{\varphi/h} u\|^2
\end{eqnarray}
\end{proposition}
\noindent\textbf{Proof}.  
We start by modifying the weight as follows:
if $\varphi_0=\varphi: \bar{M} \to \R$ is a real valued harmonic Morse function with critical points $\{p_1,\dots,p_N\}$ in the interior $M$, we 
let $\varphi_j:\bar{M} \to \R$ be harmonic functions such that $p_j$ is not a critical point of $\varphi_j$ for $j = 1,\dots,N$, their existence 
is insured by  Lemma \ref{poleszeros}.
For all $\epsilon >0$ we define the convexified weight by $\varphi_{\eps} := \varphi - \frac{h}{2\epsilon}(\sum_{j = 0}^N|\varphi_j|^2)$.
\begin{lemma}
Let $\Omega$ be an open chart of $M$ and $\varphi_{\eps}:\Omega \to \R$ be as above. Then for all $u\in C^\infty_0(\Omega)$ and $h>0$ small enough, the following estimate holds:
\begin{eqnarray}
\label{dbar carleman}
\frac{C}{\epsilon}\|u\| \leq \|e^{-\varphi_{\eps}/h}\bar{\pl} e^{\varphi_{\eps}/h}u\|
\end{eqnarray}
\end{lemma}
\noindent{\bf Proof}
We use complex coordinates $z=x+iy$ in the chart where $u$ is supported and then integrate by parts so that we have 
\begin{eqnarray*}
\|e^{-\varphi_{\eps}/h}\bar{\pl}e^{\varphi_{\eps}/h}u\|^2 &=& \frac{1}{4}
\Big(\|(\partial_{x} + \frac{i\pl_y\varphi_\eps}{h})u + (i\partial_{y} + \frac{\pl_x\varphi_\eps}{h})u\|^2\Big)\\
&=& \frac{1}{4}\Big(\|(\partial_{x} + \frac{i\pl_y\varphi_\eps}{h})u\|^2 + \|(i\partial_{y} + \frac{\pl_x\varphi_\eps}{h})u\|^2 + \frac{1}{h}\cjg u\Delta\varphi_\eps ,u\cjd\Big)
\end{eqnarray*}
where $\Delta:=-(\pl_x^2+\pl_y^2)$. Then $\cjg u\Delta \varphi_\eps,u\cjd = \frac{h}{\epsilon}(|d\varphi_0|^2 + |d\varphi_1|^2 + .. + |d\varphi_N|^2)|u|^2$,
since $\varphi_j$ are harmonic, so the proof follows from the fact that $|d\varphi_0|^2+|d\varphi_1|^2 + .. + |d\varphi_N|^2$ is uniformly 
bounded away from zero.\qed\\
%For the proof we will utilize the following Carleman estimate for the $\bar\partial$ operator on an open subset of $\R^2$. This was proved by Mikko Salo
%\begin{proposition}
%Let $M\subset \R^2$ be a bounded open subset. Let $\varphi: M\to\R$ be a harmonic function (no restriction on crit points) and $\mu: M\to\R$ satisfy $\Delta \mu \geq \delta >0$. Then for all $0< h << \epsilon << \delta$, we have the estimate
%\begin{eqnarray}
%\label{dbar carleman}
%\frac{c}{\epsilon}\|u\|^2 \leq \|e^{-\varphi_{\eps}/h} \bar\partial e^{\varphi_{\eps}/h} u\|^2
%\end{eqnarray}
%here $\varphi_{\eps} = \varphi + \frac{h}{\epsilon}\mu$. The constant $c$ depends only on $\varphi$ and $\mu$.
%\end{proposition}

The main step to go from (\ref{dbar carleman}) to (\ref{carleman}) is the following lemma which is a slight modification of the proof in \cite{IUY}:
\begin{lemma}\label{ce on patch}
With the same assumption as Proposition \ref{carlemanestimate} and if $\Omega$ is a chart of $(M,g)$ 
chosen sufficiently small and containing at most one critical point of $\varphi$, then we have
\begin{equation}\label{estimateDelta}
\frac{c}{\epsilon}(\frac{1}{h}\|u\|^2 + \frac{1}{h^2}\|u |d\varphi|\|^2 +  \frac{1}{h^2}\|u |d\varphi_{\eps}|\|^2+ \|du\|^2) \leq C\|e^{-\varphi_{\eps}/h}\Delta_g e^{\varphi_{\eps}/h} u\|^2
\end{equation}
or equivalently,
\[\frac{c}{\epsilon}\Big(\frac{1}{h}\|e^{-\varphi_{\eps}/h}u\|^2 + \frac{1}{h^2}\|e^{-\varphi_{\eps}/h}u |d\varphi|\|^2 +  \frac{1}{h^2}\|e^{-\varphi_{\eps}/h}u |d\varphi_{\eps}|\|^2+ \|e^{-\varphi_{\eps}/h}du\|^2\Big) \leq C\|e^{-\varphi_{\eps}/h}\Delta_g  u\|^2\]
for all $0<h \ll\epsilon \ll 1$ and $u\in C^\infty_0(M)$.
\end{lemma}
\noindent{\bf Proof}. Since the norms induced by the metric $g$ in the chart are conformal 
to Euclidean norms and $\Delta_g=-e^{2f}(\pl_x^2+\pl_y^2)=e^{2f}\Delta$ in the complex coordinate chart for some smooth function 
$f$, it suffices to get the estimate \eqref{estimateDelta} for Euclidean norms and Laplacian.
Clearly, we can assume $u\in C^\infty_0(M)$ to be real valued without loss of generality. By (\ref{dbar carleman}) we have
\[\|e^{-\varphi_{\eps}/h}\Delta e^{\varphi_{\eps}/h} u\|^2  = 4\|e^{-\varphi_{\eps}/h}\bar\partial e^{\varphi_{\eps}/h}e^{-\varphi_{\eps}/h}\partial e^{\varphi_{\eps}/h} u\|^2 \geq \frac{c}{\epsilon}\|e^{-\varphi_{\eps}/h}\partial e^{\varphi_{\eps}/h} u\|^2 =\frac{c}{\epsilon}\|\partial u + \frac{\pl\varphi_{\eps}}{h} u\|^2. \]
Using the fact that $u$ is real valued, we get that
\[\|e^{-\varphi_\eps/h}\Delta e^{\varphi_\eps/h} u\|^2 
\geq \frac{c}{\eps}\Big(\|du\|^2 + \frac{1}{h^2} \|u |d\varphi_{\eps}|\|^2 + \frac{2}{h}\cjg\pl_xu, u\pl_x\varphi_\eps\cjd+ 
\frac{2}{h} \cjg\pl_y u ,u\pl_y \varphi_{\eps}\cjd\Big)\]
Using the fact that $u$ is real valued, that $\varphi$ is harmonic and that $\sum_{j=0}^N|d\varphi_j|^2$ is uniformly 
bounded below, we see that
\begin{equation}\label{ineq1}
\frac{2}{h}\cjg\pl_xu, u\pl_x\varphi_\eps\cjd+ \frac{2}{h} \cjg\pl_y u,u\pl_y \varphi_\eps\cjd = 
\frac{1}{h}\cjg u,u \Delta\varphi_\eps\cjd \geq \frac{C}{\epsilon}\|u\|^2
\end{equation}
for some $C>0$ and therefore,
\[\|e^{-\varphi_{\eps}/h}\Delta_e e^{\varphi_{\eps}/h} u\|^2 \geq \frac{c}{\epsilon}(\|du\|^2 + \frac{1}{h^2} \|u |d\varphi_{\eps}|\|^2  + \frac{C}{\epsilon}\|u\|^2).\]
If the diameter of the chart $\Omega$ is chosen small (with size depending only on $|{\rm Hess}\varphi_0|(p)$) 
with a unique critical point $p$ of $\varphi_0$ inside, 
one can use integration by part and the fact that the critical point is non-degenerate to obtain 
\begin{equation}\label{ineq2} 
\|\bar{\pl} u\|^2 + \frac{1}{h^2} \|u |\pl \varphi_0|\|^2 \geq \frac{1}{h}\left|\int \pl_{\bar{z}}(u^2)\bbar{\pl_z\varphi_0} dxdy \right|
\geq \frac{1}{h}\left|\int u^2\,\bbar{\pl_z^2 \varphi_0}\, dxdy \right|\geq
\frac{C'}{h}\|u\|^2
\end{equation}
for some $C'>0$. Clearly the same estimate holds trivially if $\Omega$ does not contain critical point of $\varphi_0$.  
Thus, combining with \eqref{ineq1}, there are positive constants $c,c',C''$ such that for $h$ small enough
\[\|e^{-\varphi_{\eps}/h}\Delta e^{\varphi_{\eps}/h} u\|^2 \geq \frac{c}{\epsilon}(\|du\|^2 + \frac{1}{h^2} \|u |d\varphi_0|\|^2-\frac{C''}{\eps^2}\|u\|^2) \geq\frac{c'}{\epsilon}(\|du\|^2 + \frac{1}{h^2} \|u |d\varphi_0|\|^2 + \frac{1}{h}\|u\|^2 ) .\]
Combining the two above inequalities gives the desired estimate.
\qed\\

%\begin{coro}
%Let $(M,g)$ be a Riemann surface and $u\in C^\infty_0(M)$ be supported in a single conformal coordinate chart. Then,
%\[\frac{c}{\epsilon}(\frac{1}{h}\|u\|^2 + \frac{1}{h^2}\|u |d\varphi_0|\|^2 +  \frac{1}{h^2}\|u |d\varphi_{\eps}|\|^2+ \|du\|^2) \leq C\|e^{-\varphi_{\eps}/h}\Delta_g e^{\varphi_{\eps}/h} u\|^2\]
%or equivalently,
%\[\frac{c}{\epsilon}(\frac{1}{h}\|e^{-\varphi_{\eps}/h}u\|^2 + \frac{1}{h^2}\|e^{-\varphi_{\eps}/h}u |d\varphi_0|\|^2 +  \frac{1}{h^2}\|e^{-\varphi_{\eps}/h}u |d\varphi_{\eps}|\|^2+ \|e^{-\varphi_{\eps}/h}du\|^2) \leq C\|e^{-\varphi_{\eps}/h}\Delta_g  u\|^2\]
%for all $0<h <<\epsilon << 1$ 
%\end{coro}
\noindent{\bf Proof of Proposition \ref{carlemanestimate}}. Using triangular inequality and absorbing the term $||Vu||^2$
into the left hand side of \eqref{carleman}, it suffices to prove \eqref{carleman} with  $\Delta_g$ instead of $\Delta_g+V$.
Let $v \in C_0^\infty(M)$, we have by Lemma \ref{ce on patch} that there exist constants $c,c',C,C'>0$ such that
\[\begin{gathered}
\frac{c}{\epsilon}(\frac{1}{h}\|e^{-\varphi_{\eps}/h}v\|^2 + \frac{1}{h^2}\|e^{-\varphi_{\eps}/h}v |d\varphi|\|^2 +  \frac{1}{h^2}\|e^{-\varphi_{\eps}/h}v |d\varphi_{\eps}|\|^2+ \|e^{-\varphi_{\eps}/h}dv\|^2)  \leq\\
 \sum_j \frac{c'}{\epsilon}(\frac{1}{h}\|e^{-\varphi_{\eps}/h}\chi_jv\|^2 + \frac{1}{h^2}\|e^{-\varphi_{\eps}/h}\chi_jv |d\varphi|\|^2 +  \frac{1}{h^2}\|e^{-\varphi_{\eps}/h}\chi_jv |d\varphi_{\eps}|\|^2+ \|e^{-\varphi_{\eps}/h}d(\chi_jv)\|^2) \leq\\ 
 \sum_j C\|e^{-\varphi_{\eps}/h}\Delta_g  (\chi_jv)\|^2\leq  C'\|e^{-\varphi_{\eps}/h}\Delta_g  v\|^2 + C' \|e^{-\varphi_{\eps}/h} v\|^2 + C' \|e^{-\varphi_{\eps}/h} dv\|^2 
\end{gathered}
\]
where $(\chi_j)_j$ is a partition of unity associated to the complex charts on $M$. 
Since constants on both sides are independent of $\eps$ and $h$, we can take $\eps$ small enough so that 
$C' \|e^{-\varphi_\eps/h} v\|^2 + C' \|e^{-\varphi_\eps/h} dv\|^2$ can be absorbed to the left side. 
Now set $v = e^{\varphi_\eps/h} w$, we have 
\[\frac{1}{h}\|w\|^2 + \frac{1}{h^2}\|w |d\varphi|\|^2 +\frac{1}{h^2}\|w |d\varphi_{\eps}|\|^2 +\|dw\|^2 \leq C\|e^{-\varphi_{\eps}/h}\Delta_g e^{\varphi_{\eps}/h} w\|^2\]
Finally, fix $\epsilon >0$ and set $u = e^{\frac{1}{\epsilon}\sum_{j = 0}^N|\varphi_j|^2}w$ and use the fact that $e^{\frac{1}{\epsilon}\sum_{j = 0}^N|\varphi_j|^2}$ is independent of $h$ and bounded uniformly away from zero 
and above, we then obtain the desired estimate for $h \ll\epsilon$.\qed

\end{section}

\begin{section}{Complex Geometric Optics on a Riemann Surface}\label{CGOriemann}
As explained in the Introduction, the method for recovering the potential at a point $p$ is to construct complex geometric optic solutions depending on a small parameter $h>0$, with phase a Carleman weight (here a Morse holomorphic function), and such that the phase has a non-degenerate 
critical point at $p$, in order to apply the stationary phase method. 
 
First consider any continuous extension of $V$ to $M$, still denoted $V$ for simplicity. 
Choose $p\in M_0$ such that there exists a Morse holomorphic function $\Phi=\varphi+i\psi$ on $\mc{O}$, $C^k$ in $\bar{\mc{O}}$, 
with a critical point at $p$ and where $\mc{O}$ is chosen like in first section, ie. such that $M_0\subset\mc{O}\subset M$. Obviously $\Phi$ has isolated critical points in $\mc{O}$ and thus by reducing slightly $\mc{O}$ if necessary, we can assume that $\Phi$
has no critical point on its boundary $\partial\bar{\mc{O}}$. 
The purpose of this section is to construct solutions $u$ on $\mc{O}$ of $(\Delta +V)u = 0$ of the form
\begin{equation}
\label{cgo}
u = e^{\Phi/h}(a + r_1 + r_2)
\end{equation}
for $h>0$ small, where $a$ is a holomorphic function on $\mc{O}$ such that $a(p)\not=0$ and $r_1,r_2$ will be reminder terms which
are small as $h\to 0$ and have particular properties near the critical points of $\Phi$.
More precisely, $r_2$ will be a $O_{L^2}(h^{3/2-\eps})$ for all $\eps>0$ and $r_1$ will be a $O_{L^2}(h^{1-\eps})$ but with 
an explicit expression, which can be used to obtain sufficient informations from the stationary phase method.

\begin{subsection}{Construction of $r_1$}\label{constr1}
For all $\eps>0$ we want to construct $r_1$ which satisfies
\[e^{-\Phi/h}(\Delta_g +V)e^{\Phi/h}(a + r_1) = O(h^{1-\eps})\]
in $L^2$ and $\|r_1\|_{L^2} = O(h^{1-\eps})$.
We let $G$ be the Green operator of the Laplacian on the smooth surface with boundary 
$\bar{\mc{O}}$ with Dirichlet condition, so that $\Delta_gG={\rm Id}$ on $L^2(\mc{O})$. In particular this implies 
that $\bar{\partial}\partial G=\frac{i}{2}\star^{-1}$ where $\star^{-1}$ is the inverse of $\star$ mapping functions to $2$-forms.
First, we will search for $r_{1}$ satisfying
\begin{equation}
\label{dequation}
e^{-2i\psi/h}\partial e^{2i\psi/h} r_1 = -\pl G (aV) + \omega + O(h^{1-\eps})
\end{equation}
in $H^1(\mc{O})$, with $\omega$ a holomorphic 1-form on $\mc{O}$  and $\|r_1\|_{L^2} = O(h^{1-\eps})$. 
Indeed, using the fact that $\Phi$ is holomorphic we have
\[e^{-\Phi/h}\Delta_ge^{\Phi /h}=-2i\star  \bar{\pl} e^{-\Phi/h}\pl e^{\Phi/h}=-2i\star  \bar{\pl} e^{-\frac{1}{h}(\Phi-\bar{\Phi})}\pl e^{\frac{1}{h}(\Phi-\bar{\Phi})}=
-2i\star \bar{\pl}e^{-2i\psi/h}\pl e^{2i\psi/h}\]
and applying $-2i\star\bar{\pl}$ to \eqref{dequation}, this gives
\[e^{-\Phi/h}(\Delta_g+V)e^{\Phi/h}r_1=-aV+O_{L^2}(h^{1-\eps}).\]
The form $\omega$ above, will be chosen as a correction term to optimize the use of the stationary phase later, this is why we need the following
\begin{lemma}\label{formomega}
Let $\{p_0,...,p_N\}$ be finitely many points on $\mc{O}$ and let $g$ be a continuous section of $T^*_{1,0}\mc{O}$.
Then there exists a holomorphic $1$-form $\omega$ on $\mc{O}$ such that $(g-\omega)(p_i)=0$ for all $i=0,\dots,N$.
\end{lemma}
\noindent{\bf Proof}. First by Riemann-Roch formula \eqref{riemannroch}, there exists a meromorphic $1$-form $v$ on 
$M$, holomorphic on $\mc{O}$, which has a zero of order greater or equal to $1$ at all $p_1,\dots,p_N$, so using Lemma 
\ref{poleszeros}, we can multiply it by a meromorphic function $f_j$ on $\mc{O}$, holomorphic on $\mc{O}\setminus\{p_j\}$, 
with a pole of order exactly $n_j$ at $p_j$ if $n_j$ is the order of $p_j$ as a zero of $v$,
so that $v_j:=f_jv$ is a holomorphic $1$-form on $\mc{O}$ with no zero at $p_j$ and zeros of order larger or equal to $1$
at all other $p_k$ for $k\not=j$. Now since $T^*_{1,0}\mc{O}$ is a  complex line bundle, 
there is a complex number $c_j\in\cc$ such that $g(p_j)=c_jv_j(p_j)$. Thus it is clear that
$\omega=\sum_{j=1}^Nc_jv_j$ satisfies the claim.\qed\\    

With this lemma, we choose $\omega$ to be a holomorphic $1$-form on $\mc{O}$ such that at all 
critical point $p'$ of $\Phi$, we have $(\pl G(aV) - \omega)(p') = 0$; this can be done since
$\pl G(aV)$ is a $1$-form with value in $T^*_{1,0}\mc{O}$.\\
 
We will construct $r_1= r_{1,1} + r_{1,2}$ in two steps. 
First, we will construct $r_{1,1}$ to solve equation \eqref{dequation} 
locally near critical points of $\Phi$ by using coordinate charts. Then we will 
construct the global correction term $r_{1,2}$ away from critical points.\\

Let $p'$ be a critical point of $\Phi$ and ${\cal U}(p')$ be a complex coordinate chart $z$ 
containing $p'$ but no other critical points of $\Phi$. 
In local coordinates one has $-\pl G (aV) + \omega = b(z) dz$ 
for some $C^{3}$ function $b$ vanishing at $p'$. Let $\chi_1\in C_0^\infty({\cal U}(p'))$ 
such that $\chi_1=1$ in a neighbourhood of $p'$ 
and let $\chi\in C^{\infty}_0({\cal U}(p'))$ with $\chi= 1$ on an open set containing the support of $\chi_1$. 
Define for $z\in {\cal U}(p')$
\begin{equation}\label{defr11}
r_{1,1} (z):= \chi(z)e^{-2i\psi/h}R(e^{2i\psi/h}\chi_1b)(z)
\end{equation}
where $Rf(z) := \int_{\R^2} \frac{1}{z-\xi}f d\xi\wedge d\overline\xi$ for $f\in L^\infty$ compactly supported 
is the classical Cauchy-Riemann operator. Extend $r_{1,1}$ trivially outside of ${\cal U}(p')$. Then
\begin{equation}\begin{gathered}\label{r11}
e^{-2i\psi/h}\pl(e^{2i\psi/h}r_{1,1}) = \chi_1(-\pl G(aV) + \omega) + \eta\\
\textrm{ with }\eta:= e^{-2i\psi/h}R(e^{2i\psi/h}\chi_1b)\wedge \pl\chi \end{gathered}
\end{equation}
where the form $\eta$ makes sense globally on $\mc{O}$ since $\pl\chi$ is supported in ${\cal U}(p')$. 
Note that $\eta$ is a $C^4$ form with value in $T^*_{1,0}\mc{O}$. 
Now the support of $\pl\chi$, thus of $\eta$, is contained in the complement of the support of $\chi_1$. 
By stationary phase and the fact that $b = 0$ at all critical points of $\Phi$, one has (combine for instance 
Propositions 3.2 and 3.4 of \cite{IUY}) 
\begin{equation}
\label{asymptoticeta}
\|\eta\|_{\infty} \leq Ch^2\ \ \ 
\textrm{ and } \ \ \ \|\Delta \eta\|_{\infty}\leq C.
\end{equation}
The term $r_{1,1}$ is supported in ${\cal U}(p')$ for a fixed critical point $p'$ and depends 
on $p'$, let us write it $r^{p'}_{1,1}$ instead, but since our discussion 
did not depend on the choice of $p'$, we can sum the $r_{1,1}^{p'}$ 
over the critical points $p'$ to define a term $r_{1,1}$.\\ 

Next we define $r_{1,2}$ by the equation
\[r_{1,2}\pl\Phi = h\Big(-\eta + (1-\chi_1)(-\pl G(aV)+\omega)\Big).\]
so that 
\begin{equation}\label{r12}
e^{-\Phi/h}\pl e^{\Phi/h}r_{1,2}=\pl r_{1,2}-\eta+(1-\chi_1)(-\pl G(aV)+\omega).
\end{equation}
There is a well defined $C^{3}$ function $r_{1,2}$ satisfying this equation 
since both $\pl\Phi$ and the right hand side have values in the bundle 
$T^*_{1,0}\mc{O}$ and moreover the right hand side has support which does not intersect the critical points of $\Phi$.

We now derive the asymptotic properties of $r_{1,1}$ and $r_{1,2}$
\begin{lemma}
\label{asymptoticsr11r12}
For all $\eps >0$, the following estimates hold
\[\begin{gathered}
\|r_{1,1}\|_{L^2}\leq Ch^{1-\eps}, \ \ \
\|r_{1,2}\|_{\infty} \leq Ch^1\ \ \ \textrm{ and } \ \ \ \|\Delta r_{1,2}\|_{\infty}\leq Ch.
\end{gathered}\]
\end{lemma}
\noindent{\bf Proof}. The first estimate is a local result and comes from classical properties of $R$ 
as proved in Proposition 3.5 of \cite{IUY}. The ones involving $r_{1,2}$ follow directly from (\ref{asymptoticeta}).
\qed
\begin{lemma}
\label{solve equation to next order}
%\[\Delta r^{1} + \frac{1}{h}*\overline\partial r^{1}\wedge\partial\Phi = aV + R\]
%where $\|R\|_\infty\leq Ch$.
With $r_1:=r_{1,1}+r_{1,2}$ constructed above, then for all $\eps>0$
\[e^{-\Phi/h}(\Delta +V)e^{\Phi/h}(a + r_1) = O_{L^2}(h^{1-\eps}).\]
\end{lemma}
\noindent{\bf Proof}. First, we write
\[e^{-\Phi/h}\pl e^{\Phi/h}r_1= e^{-2i\psi/h}\pl(e^{2i\psi/h}r_{1,1})+ (\pl + \frac{1}{h}\partial\Phi)r_{1,2}.\]
and by \eqref{r11} and \eqref{r12} this implies
\[e^{-\Phi/h}\pl e^{\Phi/h}r_1=-\pl G(aV)+\omega+\pl r_{1,2}\]
and applying $-2i\star\bar{\pl}$
\[e^{-\Phi/h}\Delta e^{\Phi/h}r_1=-aV+\Delta r_{1,2}.\] 
The proof is complete since we know from Lemma \ref{asymptoticsr11r12} 
that $\|\Delta r_{1,2}\|_\infty \leq Ch$.\qed
\end{subsection}

\begin{subsection}{Construction of $r_2$}
The goal of this section is to complete the construction of the complex geometric optic solutions by the following proposition:
\begin{proposition}
\label{completecgo}
For all $\eps >0$ there exist solutions to $(\Delta +V)u = 0$ of the form (\ref{cgo}) with $r_1=r_{1,1}+r_{1,2}$ 
constructed in the previous section and $r_2$ satisfying $\|r_2\|_{L^2}\leq Ch^{3/2-\eps}$
\end{proposition}
This is a consequence of the following Lemma (which follows from the Carleman estimate obtained above) 
\begin{lemma}
\label{standardargument}
Let $q\in L^\infty(\mc{O})$ and $f\in L^2(\mc{O})$. For all $h>0$ small enough, there exists a solution $v\in L^2$ to the equation
\[e^{-\varphi/h}(\Delta_g +V) e^{\varphi/h}v = f\]
satisfying
\[\|v\|_{L^2} \leq Ch^\demi\|f\|_{L^2}\]
\end{lemma}
\noindent{\bf Proof}. The proof is the same than Proposition 2.2 of \cite{IUY}, we repeat the argument for the convenience of the reader.
 Define for all $h >0$ the real vector space $\mc{A}:=\{u\in H_0^1(\mc{O}); (\Delta_g+V)u\in L^2(\mc{O})\}$
%${\cal A} := \{e^{\phi/h}(\Delta +V)e^{-\phi/h} u \mid u\in C^{\infty}_0(\mc{O})\}$
equipped with the real scalar product 
\[(u,w)_{\mc{A}}:=\int_{\mc{O}}e^{-2\varphi/h}(\Delta_gu+Vu)(\Delta_gw+Vw)dg.\] 
By the Carleman estimate of Proposition \eqref{carlemanestimate}, the space $\mc{A}$ is a Hilbert space equipped with the scalar product above and
so the linear functional $L:w\to \int_{\mc{O}}e^{-\varphi/h}fw \,dg$ on $\mc{A}$ is continuous and norm bounded by 
$h^\demi||f||_{L^2}$ by Proposition \eqref{carlemanestimate},  
and by Riesz theorem there is an element $u\in\mc{A}$ such that $(.,u)_{\mc{A}}=L$  
and with norm bounded by the norm of 
$L$. It remains to take $v:=e^{-\varphi/h}(\Delta_g u+Vu)$ which solves $\Delta_gv+Vv=f$  and 
which in addition satisfies the desired norm estimate.\qed\\

\noindent{\bf Proof of Proposition \ref{completecgo}}. We note that $(\Delta +V)e^{\Phi/h}(a + r_1 + r_2) = 0$ if and only if
\[e^{-\Phi/h}(\Delta +V)e^{\Phi/h}r_2 = - e^{-\Phi/h}(\Delta +V)e^{\Phi/h}(a + r_1)\]
By Lemma \ref{standardargument} one can find such an $r_2$ which satisfies
\[\|r_2\|_{L^2} \leq Ch^\demi\|e^{-\Phi/h}(\Delta +V)e^{\Phi/h}(a + r_1)\|_{L^2} \leq Ch^{3/2 - \eps}\]
where the last inequality comes from Lemma \ref{solve equation to next order}.\qed

\end{subsection}
\end{section}
\begin{section}{Recovering the potential}
We now assume that $V_1,V_2\in C^{2}(\bar{M_0})$ are two real valued potentials such that the respective Caucy data spaces
$\mc{C}_1,\mc{C}_2$ for the operators $\Delta_g+V_1$ and $\Delta_g+V_2$ are equal.
Let $p\in M_0$ and $\mc{O}$ with $M_0\subset \mc{O}\subset M\setminus\{q\}$ such that, using Proposition \ref{criticalpoints}, we can choose a holomorphic Morse function $\Phi=\varphi+i\psi$ on $\mc{O}$, $C^k$ in $\bar{\mc{O}}$ for some large $k\in\nn$, with a critical point at $p$. 
By reducing slightly $\mc{O}$ if necessary, we can assume that $\Phi$ has no critical points on $\pl\bar{\mc{O}}$ and finitely many critical points 
in $\mc{O}$.
\begin{proposition}
\label{identcritpts}
If the Cauchy data spaces agree, i.e. if $\mc{C}_1 = \mc{C}_2$, then $V_1(p)= V_2(p)$.
\end{proposition}
\noindent{\bf Proof}. By boundary identifiability (see for example \cite{DKSU}), 
one has $V_1 = V_2$ on $\pl M_0$ to second order and therefore we can extend $V_1$, 
$V_2$ to be $C^2$ to $\bar{\mc{O}}$ such that they agree outside of $\pl M_0$.\\
Let $a$ be a holomorphic function on $\mc{O}$ with $a(p)\neq 0$ and $a(p') = 0$ for all other critical point $p'$ of $\Phi$. 
The existence is insured by Lemma \ref{poleszeros} as follows: by Riemann-Roch, we can find a holomorphic function on $M\setminus\{q\}$ such that $a(p') = 0$ for all $p'\neq p$. Either at $p$ this function does not vanish and we have our function $a$, or there is a zero of 
order say $N$, in which case one can multiply it by a meromorphic function on $M$, holomorphic on $M\setminus\{q,p\}$ 
with a pole of order exactly $N$ at $p$ (the existence of which is proved in Lemma \ref{poleszeros}).  
Let $u_1$ and $u_2$ be $H^2$ solutions on $\bar{\mc{O}}$ to 
\[(\Delta_g +V_j)u_j = 0\]
constructed in Section \ref{CGOriemann} with $\Phi$ for Carleman weight for $u_1$ and $-\Phi$ for $u_2$, thus of the form
\[u_1 = e^{\Phi/h}(a + r_1^{1} + r_2^{1}), \quad u_2=e^{-\Phi/h}(a+r_1^{2}+r_2^{2})\]
and with boundary value $u_j|_{\pl M_0}=f_j$.
Since $\bbar{u_2}$ is also a solution, we can write by Green formula
\[\begin{split}
\int_{M_0}u_1(V_1 - V_2) \bbar{u_2} {\rm dv}_g&=- \int_{M_0} (\Delta_gu_1. \bbar{u_2} - u_1.\Delta_g\bbar{u_2}){\rm dv}_g\\
&=-\int_{\pl M_0} (\pl_n u_1.\bbar{f_2}-f_1. \pl_n \bbar{u_2}) {\rm dv}_{g}
 \end{split}\]
Since the Cauchy data for $\Delta_g + V_1$ agrees with that of $\Delta_g + V_2$, there exists a solution $v$ of the boundary value problem
\[(\Delta_g + V_2) v = 0,\ \ \ \ \ \ \ v|_{\pl M_0}= f_1\]
satisfying $\pl_n v = \pl_n u_1$ on $\partial M_0$. This implies that 
\[\begin{split}
\int_{M_0}u_1(V_1 - V_2) \bbar{u_2} {\rm dv}_g&=- \int_{M_0} (\Delta_gu_1. \bbar{u_2} - u_1\Delta_g\bbar{u_2}){\rm dv}_g=-\int_{\pl M_0} (\pl_n u_1\bbar{f_2}-f_1 \pl_n \bbar{u_2}) {\rm dv}_{g}\\
&=-\int_{\pl M_0} (\pl_n v \bbar{f_2} - v\pl_n \bbar{u_2}){\rm dv}_{g}= -\int_{M_0} (\Delta_gv. \bbar{u_2} - v\Delta_g\bbar{u_2}){\rm dv}_g = 0\\
 \end{split} \]
since $\Delta_g + V_2$ annihilates both $v$ and $u_2$. Then by using the estimates in Lemma \ref{asymptoticsr11r12} and Proposition \ref{completecgo} that we have, as $h\to 0$,
\[\int_{M_0}e^{2i\psi/h}|a|^2(V_1 - V_2) {\rm dv}_g +  
\int_{M_0}e^{2i\psi/h}(\bbar{a} r_1^{1} + a \bbar{r_1^{2}})(V_1 - V_2){dv}_g + O(h^{3/2 -\eps}) = 0\]
for all $\eps>0$.
By stationary phase the first term can be developed as follows  
\[\int_{M_0}e^{2i\psi/h}|a|^2(V_1 - V_2) {\rm dv}_g = 
\int_{\mc{O}}e^{2i\psi/h}|a|^2(V_1 - V_2){\rm dv}_g = Ch(V_1(p) - V_2(p)) + O(h^2)\]
for some $C\neq 0$.  Therefore,
\[ Ch(V_1(p) - V_2(p))+  \int_{\mc{O}}e^{2i\psi/h}(\bbar{a} r_1^{1} + a \bbar{r_1^{2}})(V_1 - V_2){\rm dv}_g + O(h^{3/2 -\eps}) = 0.\]
It suffices then to show that 
\[\int_{M_0}e^{2i\psi/h}(\bbar{a} r_1^{1} + a \bbar{r_1^{2}})(V_1 - V_2){\rm dv}_g = o(h).\]
This can be accomplished by the following argument: let us deal with the $\bar{a}r_{1}^1$ term since the other one is 
exactly similar, then we localize near critical points and by stationary phase,
\[\int_{\mc{O}}e^{2i\psi/h}\bar{a} r_1^{1}(V_1 - V_2){\rm dv}_g = 
\sum_{p'\in \mc{O}, d\Phi(p')=0}
\int_{{\cal U}(p')}e^{2i\psi/h}\chi_{p'} \bar{a}r_1^{1}(V_1 - V_2){\rm dv}_g + O(h^2)\]
where ${\cal U}(p')$ is a small neighbourhood of $p'$ that we used to define $r_{1,1}$ in the Subsection \ref{constr1} and $\chi_{p'}$
the smooth cutoff functions supported in that neighbourhood. First we observe that by stationary phase and the fact that $\|r^{1}_{1,2}\|_\infty \leq Ch$ (Lemma \ref{asymptoticsr11r12}), we have
\begin{eqnarray*}
\int_{{\cal U}(p')}e^{2i\psi/h}\chi_{p'} \bar{a} r_1^{1}(V_1 - V_2){\rm dv}_g &=&\int_{{\cal U}(p')}e^{2i\psi/h}\chi_{p'} \bar{a} (r_{1,1}^{1} + r_{1,2}^{1})(V_1 - V_2){\rm dv}_g\\ %+ \int_{{\cal U}(p')}e^{i\psi/h}\chi_{p'} \bar{a} r_{1,2}^{1}(V_1 - V_2){\rm dv}_g\\
&=&\int_{{\cal U}(p')}e^{2i\psi/h}\chi_{p'} \bar{a} r_{1,1}^{1}(V_1 - V_2){\rm dv}_g + O(h^2)
\end{eqnarray*}

 The complex coordinates are denoted $z=x+iy$ in these charts
and the volume form is written $\sqrt{\det(g)(z)}dxdy$.   
For each critical points $p'$ using the local representation in ${\cal U}(p')$ of $r_1^{1}$ in \eqref{defr11} (and 
the same notations for $R$ and $b$), 
we obtain
\begin{eqnarray*}
\int_{{\cal U}(p')}e^{i\psi/h}\chi_{p'} \bar{a} r_1^{1}(V_1 - V_2){\rm dv}_g &=&\int_{{\cal U}(p')}e^{2i\psi/h}\chi_{p'} \bar{a} r_{1,1}^{1}(V_1 - V_2){\rm dv}_g + O(h^2)\\ &=& \int_{\R^2}
\chi_{p'}^2 R(e^{2i\psi/h}\chi_1b)a(V_1 - V_2) \sqrt{\det g}\, dxdy + O(h^2)\\ &=& 
\int_{\R^2}e^{2i\psi/h}\chi_1 bR^*\Big(\chi_{p'}^2 a(V_1 - V_2)\sqrt{\det g}\Big)\,dxdy + O(h^2).
\end{eqnarray*}
where $R^*$ is the adjoint of $R$ on $L^2$ compactly supported functions in $\rr^2$, 
which has the same mapping properties than $R$ on $C^k$ (its integral kernel 
being of the same form).
Since $b$ and $R^*(\chi_{p'}^2 a(V_1 - V_2)\sqrt{\det g})$ are $C^3$ in \eqref{defr11} and $b$ vanishes at $p'$, 
we obtain by stationary phase that the above integral is $O(h^2)$, which completes the proof.\qed\\
\noindent{\bf Proof of Theorem \ref{identif}}.
By combining Proposition \ref{identcritpts} and the Proposition \ref{criticalpoints}, we show that $V_1=V_2$
on a dense set of $M_0$, so by continuity of $V_1,V_2$ they agree everywhere.\qed
\end{section}

\section{Inverse scattering}\label{invscat}

We first obtain, as a trivial consequence of Theorem \ref{identif}, a result about inverse scattering for 
asymptotically hyperbolic surface.

\noindent{\bf Proof of Corollary \ref{coroAH}}. Let $x$ be a smooth boundary defining function of $\pl\bar{X}$, and let $\bar{g}=x^2g$
be the compactified metric and define $\bar{V}_j:=V_j/x^2 \in C^{\infty}(\bar{X})$. 
By conformal invariance of the Laplacian in dimension $2$, one has 
\[\Delta_g+V_j=x^2(\Delta_{\bar{g}}+\bar{V}_j)\]
and so if $\ker_{L^2}(\Delta_g+V_1)=\ker_{L^2}(\Delta_g+V_2)$ and $\mc{S}_1=\mc{S}_2$, then 
the Cauchy data space 
\[\{(u|_{\pl\bar{X}},\pl_xu|_{\pl\bar{X}})\in C^{\infty}(\pl\bar{X})\x C^{\infty}(\pl\bar{X}); 
u \in C^{\infty}(\bar{X}), (\Delta_{\bar g}+\bar{V}_j)u=0\}\]
for the operators $\Delta_{\bar{g}}+\bar{V}_j$ are obviously the same. Then it suffices to apply the result in Theorem \ref{identif}.
\qed\\

Now we consider the asymptotically Euclidean scattering at $0$ frequency. Recall that an asymptotically Euclidean surface $(X,g)$
has a smooth compactification $\bar{X}$ with $g$ of the form 
\[g=\frac{dx^2}{x^4}+\frac{h(x)}{x^2}\]
near the boundary $\pl\bar{X}$, where $x$ is a smooth boundary defining function and 
$h(t)$ a smooth family of metrics on the slice $\{x=t\}\simeq \pl\bar{X}$ down to $\{x=0\}=\pl\bar{X}$ and such that $h(0)=d\theta_{S^1}^2$ 
is the metric with length $2\pi$ on each of the, say $L$, copies of $S^1$ composing the boundary.
In particular it is conformal to an asymptotically cylindrical metric, or 'b-metric' in the sense of Melrose \cite{APS}, 
\[g_b:=x^2g=\frac{dx^2}{x^2}+h(x)\]
and the Laplacian satisfies $\Delta_g=x^2\Delta_{g_b}$. Each end of $X$ is of the form $(0,\eps)_x\x S^1_{\theta}$ 
and the operator $\Delta_{g_b}$ has the expression in the ends
\[\Delta_{g_b}=-(x\pl_x)^2+\Delta_{\pl\bar{X}}+ xP(x,\theta;x\pl_x,\pl_\theta)\]
for some smooth differential operator $P(x,\theta;x\partial_x, \partial_\theta)$ in the vector fields $x\pl_x, \pl_\theta$  
down to $x=0$. 
Let us define $V_b:=x^{-2}V$, which is compactly supported and 
\[H^{2m}_b:=\{u\in L^2(X,{\rm dvol}_{g_b}); \Delta^m_{g_b}u\in L^2(X,{\rm dvol}_{g_b})\}, \quad m\in\nn_0.\] 
We also define the following spaces for $\alpha\in\rr$ 
\[F_\alpha:=\ker_{x^\alpha H^2_b}(\Delta_{g_b}+V_b).\]
Since the eigenvalues of $\Delta_{S^1}$ are $\{j^2; j\in\nn_0\}$, 
the relative Index theorem of Melrose \cite[Section 6.2]{APS} shows that $\Delta_{g_b}+V_b$ 
is Fredholm from $x^{\alpha}H^2_b$ to $x^{\alpha}H^0_b$
if $\alpha\notin \zz$.
Moreover, from subsection 2.2.4 of \cite{APS}, we have that any solution of $(\Delta_{g_b}+V_b)u=0$
in $x^{\alpha}H^2_b$ has an asymptotic expansion of the form 
\[u\sim \sum_{j>\alpha, j\in \zz}\sum_{\ell=0}^{\ell_j}x^{j}(\log x)^\ell u_{j,\ell}(\theta), \quad \textrm{ as }x\to 0\]
for some sequence $(\ell_j)_j$ of non negative integers and some smooth function $u_{j,\ell}$ on $S^1$. 
%\[\textrm{with }u_j\in  E_j,  \textrm{ and }E_j:=\ker(\Delta_{\pl\bar{X}}-j^2).\]
In particular, it is easy to check that $\ker_{L^2(X,{\rm dvol}_g)}(\Delta_g+V)=F_{1+\eps}$ for $\eps\in (0,1)$.

\begin{theorem} 
Let $(X,g)$ be an asymptotically Euclidean surface and $V_1,V_2$ be two compactly supported smooth potentials
and $x$ be a boundary defining function. Let $\eps\in(0,1)$
and assume that for any $j\in\zz$ and any function $\psi\in \ker_{x^{j-\eps}H^2_b}(\Delta_{g}+V_1)$ there is a 
$\varphi\in\ker_{x^{j-\eps}H^{2}_b}(\Delta_g+V_2)$ such that $\psi-\varphi=O(x^\infty)$, and conversely. 
Then $V_1=V_2$. 
\end{theorem}
\noindent\textbf{Proof}. The idea is to reduce the problem to the compact case. 
First we notice that by unique continuation, $\psi=\varphi$ where $V_1=V_2=0$.
%With the notations
%used before, take $u_i\in \ker_{x^{-j-\eps}H^2_b}(\Delta_g+V_i)$ for $i=1,2$ and some $j\in\nn_0$.
%The solutions $u_i$ has an expansion of the form 
%\[u_i=\sum_{k=1}^ja_{-j,i}x^{-j}+\log(x)a_{-0,i}+\sum_{k=0}^j(\mc{S}_ia_{-,i})_{j}x^j+O(x^{j+1})\]
%with $a_{-,i}:=(a_{k,i})_{k}$, and $(\mc{S}_ia_{-,i})_j$ denotes the component on $E_j$ of $\mc{S}_ia_{-,i}$.
%Then if $V_1=V_2$ in $x<\eta$ for some small $\eta>0$, Green's formula in $\eps\leq x\leq \eta$ gives
%and letting $\eps\to 0$ gives directly that, with the dot product used in \eqref{symmetryS}, 
%\[ 0= \int_{\pl\bar{X}}(\mc{S}_1a_{-,1}.a_{-,2}-\mc{S}_2a_{-,2}.a_{-,1}){\rm dvol}_{\pl\bar{X}}
%+\int_{x=\eta}(x\pl_x u_1)u_2-u_1(x\pl_xu_2){\rm dvol}_{g_b}.\]
%Now if $\mc{S}_1=\mc{S}_2$ the symmetry in \eqref{symmetryS} shows that the first integral is $0$,
%and using the conformal relation $x^2g=g_b$ to rewrite the second integral in terms of ${\rm dvol}_g$, we see that 
%\[\int_{x=\eta}(\pl_n u_1)u_2-u_1(\pl_n u_2){\rm dvol}_{g}=0\]  
%where $\pl_n=x^2\pl_x$ is the normal interior pointing vector field to $x=\eta$ for the metric $g$.
Now it remains to prove that, if $R_\eta$ denote the restriction of smooth functions on $X$ to 
$\{x\geq \eta\}$ and $V$ is a smooth compactly supported potential in $\{x\geq \eta\}$,
then the set $\cup_{j=0}^\infty R_{\eta}(F_{-j-\eps})$ is dense in the set $N_V$ of $H^2(\{x\geq \eta\})$ solutions
of $(\Delta_g+V)u=0$. The proof is well known for positive frequency scattering (see for instance Lemma 3.2 in \cite{Mel}), %marker
here it is very similar so we do not give much details. The main argument is to show that it converges in $L^2$ sense
and then use elliptic regularity; the $L^2$ convergence can be shows as follows: let $f\in N_V$ such that
\[\int_{x\geq \eta}f\psi{\rm dvol}_g=0, \, \forall \, \psi \in \cup_{j=0}^\infty F_{-j-\eps},\]
then we want to show that $f=0$. By Proposition 5.64 in \cite{APS}, there exists $k\in\nn$ and a generalized right inverse $G_b$ for
$P_b=\Delta_{g_b}+V_b$ (here, as before, $x^2V_b=V$) in $x^{-k-\eps}H^2_b$, such that
$P_bG_b={\rm Id}$. This holds in $x^{-k-\eps}H^2_b$ for $k$ large enough since the cokernel of $P_b$ on this space becomes $0$ 
for $k$ large. Let $\omega=G_bf$ so that $(\Delta_{g_b}+V_b)\omega=f$, and in particular this function is $0$ 
in $\{x<\eta\}$. The asymptotic behaviour of the integral kernel $G_b(z,z')$ of $G_b$ as $z\to \infty$ 
is given in Proposition 5.64 of \cite{APS} uniformly in $z'\in\{x\geq \eta\}$, 
we have for all $J\in\nn$ and using the radial coordinates $(x,\theta)\in (0,\eps)\x S^1$ 
for $z$ in the ends
\[G_b(z,z')= \sum_{j=-k}^{J}\sum_{\ell=0}^{\ell_j} x^{j}(\log x)^{\ell}\psi_{j}(\theta,z')+o(x^{J})\]  
for some functions $\psi_{j,\ell}\in x^{k-j-\eps}H^2_b$ and some sequence $(\ell_j)_j$ of non-negative integers. 
But the fact that $(\Delta_{g_b}+V_b)G_b(z,z')=\delta(z-z')$ as distributions implies directly 
that $(\Delta_{g_b}+V_b)\psi_j(\theta,.)=0$.
Using our assumption on $f$, we deduce that $\int_{X}\psi_j(\theta,z')f(z'){\rm dvol}_{g_b}=0$ for all $j\in \nn_0$ 
and so the function $\omega$ vanish faster than all power of $x$ at infinity. Then by unique continuation, 
we deduce that $\omega=0$ in $\{x\leq \eps\}$. Since now $\omega\in H^2$, its Cauchy 
data at $x=\eta$ are $0$ and $\Delta_{g_b}+V_b$ is self adjoint for the measure ${\rm dvol}_{g_b}$, 
we can use the Green formula to obtain
\[\int_{x\geq \eta}|f|^2{\rm dvol}_{g_b}=\int_{x\geq \eta}\omega(\Delta_{g_b}+V_b)\bar{f}{\rm dvol}_{g_b}=0.\]
The $H^2$ density is easy using elliptic regularity.\qed 

\subsection{Acknowledgements}
This work started during a summer evening in Pisa thanks to the hospitality of M. Mazzucchelli and A.G. Lecuona. We thank Mikko Salo, Eleny Ionel and Rafe Mazzeo for pointing out helpful references. 
C.G. thanks MSRI and the organizers of the 'Analysis on Singular spaces' 2008 program  
for support during part of this project. This works was achieved while C.G. was 
visiting IAS under an NSF fellowship number No. DMS-0635607. 
L.T is supported by NSF Grant No. DMS-0807502.

\end{document}